\documentclass[11pt]{article}
\usepackage{amsmath}
\usepackage{amsthm,amsfonts,amssymb,epsfig,graphics,color}
\usepackage[utf8]{inputenc}
\usepackage[english]{babel}
\usepackage{hyperref}

\newcommand{\bbR}{{\mathbb{R}}}

\def\?{$^{***}$\marginpar{?}}

\newtheorem*{ques*}{Question}
\newtheorem*{prop*}{Proposition}
\newtheorem*{conj*}{Conjecture}
\newtheorem*{theo*}{Theorem}

\newtheorem{affi*}{Affirmation}

\newtheorem*{lemm*}{Lemma}
\def\?{\footnote{?}}

\newlength{\espaceavantspecialthm}
\newlength{\espaceapresspecialthm}
{\vskip \espaceapresspecialthm}

\newenvironment{defi*}{\vskip \espaceavantspecialthm \noindent \textbf{Definition.} }%
{\vskip \espaceapresspecialthm}

\newenvironment{rema}[1][]{\refstepcounter{coro} 
\vskip \espaceavantspecialthm \noindent \textbf{Remark~\thecoro
#1.} }%
{\vskip \espaceapresspecialthm}

\title{Projective transformations of rotation sets}
\author{Fran\c cois B\'eguin, Sylvain Crovisier\footnote{S.C was partially supported by the ERC project 692925 NUHGD.} and Fr\'ed\'eric Le Roux}

\begin{document}
\maketitle

\begin{abstract}
We give a new proof and extend a result of J. Kwapisz: whenever a set $C$ is realized as the rotation set of some torus homeomorphism, the image of $C$ under certain projective transformations is also realized as a rotations set.
\end{abstract}

The concept of \emph{rotation set}, introduced by M. Misiurewicz and K. Ziemian in \cite{MisZie}, is one of the most important tools to study the global dynamics of homeomorphisms of the torus $\mathbb{T}^2$. If $f$ is a homeomorphism of $\mathbb{T}^2$ isotopic to the identity, and $F$ is a lift of $f$ to $\mathbb{R}^2$, the rotation set of $F$ is a compact convex subset of the plane which describes ``at what speeds and in what directions the orbits of $f$ rotate around the torus". One of the main problems in the theory is to determine which compact convex subsets of $\mathbb{R}^2$ can be realized as the rotations sets of some torus homeomorphisms. For compact convex subsets with empty interiors (\emph{i.e.} singletons and segments), a conjectural answer to the problem has been formulated by J. Franks and M. Misiurewicz (see \cite{FraMis}). Fifteen years ago, J. Kwapicz has introduced a technical tool which allows to simplify the problem. Namely, he observed that, if a compact convex set $C\subset \mathbb{R}^2$ is realized as the rotation of a certain torus diffeomorphism, and if a projective transformation $L$ maps $C$ to a bounded set of the plane, then $L(C)$ can be realized as the rotation of another torus diffeomorphism (see \cite[section 2]{Kwa}). 

Kwapisz's proof requires to consider the suspension of the initial torus homeomorphism, and to apply a theorem of D. Fried to find a new surface of section for this flow, in the appropriate cohomology class. Fried's theorem works only for $C^1$ flows; this forces Kwapisz to consider only rotation sets of $C^1$ diffeomorphisms, whereas the natural setting for his result would be rotation sets of homeomorphisms. The purpose of the present note is to provide a more elementary proof of Kwapisz's result. Our proof remains at the level of surfaces homeomorphisms, \emph{i.e.} does not require to consider a flow on a three-dimensional manifold. It does not make use of Fried's theorem (in some sense, we replace it by the more classical fact that the only surface with fundamental group isomorphic to $\mathbb{Z}^2$ is the torus $\mathbb{T}^2$). As a consequence, it works for surface homeomorphisms without any differentiability assumption. This might be of interest in relation with some recent works related to the Franks-Misiurewicz conjecture (see~\cite{LecTal,Kor}, and the example of Avila quoted in these papers).

\begin{theo*}
Let $C$ be a compact subset of the plane which is realized as the rotation set of some torus homeomorphism. Let $L \in \mathrm{SL}(3,\mathbb{Z})$ be a projective transformation such that the image $C'$ of $C$ under $L$  is a bounded subset of the plane. Then $C'$ is also realized as the rotation set of some torus homeomorphism. 
\end{theo*}

In this statement, we use the usual affine chart to embed the plane in the projective plane. The requirement that the image of $C$ under $L$ is a bounded subset of the plane means that we demand that $L(C)$ does not meet the line at infinity. A more precise version of the above theorem will be given below.

We now recall the classical definition of the rotation set by Misiurewicz and Ziemian. We consider a self-homeomophism $f$ of the torus $\mathbb{T}^2 = \mathbb{R}^2/\mathbb{Z}^2$, and a lift $F:\mathbb{R}^2 \to \mathbb{R}^2$.  We assume that $f$ is isotopic to the identity which amounts to say that $F$ commutes with the deck transformations $S : (x,y)\mapsto (x+1,y)$ and $T : (x,y)\mapsto (x,y+1)$. The rotation set of $F$ is defined as the set of $w \in \bbR^2$ such that there exists a sequence $(z_k)_{k \geq 0}$ of points of the plane, and a sequence $(n_k)_{k \geq 0}$ of integers tending to $\infty$ such that 
$$
\frac{F^{n_k}(z_k)-z_k}{n_k}
$$
converges to $w$ as $k$ goes to infinity. 

Note that this definition depends on the choice of coordinates on the torus (in order to identify the universal cover of the torus with $\mathbb{R}^2$). In particular, it depends on the choice of a basis $(S,T)$ of the fundamental group of the torus. To make things clear we need a definition of the rotation set that makes explicit this dependence. 

\begin{defi*}\label{d.rotation-set}
Consider an action of $\mathbb{Z}^3$ on $\mathbb{R}^2$ generated by three commuting homeomorphisms $G,U,V$.
We define \emph{the rotation set of $G$ with respect to $U$ and $V$} as the set $\rho_{U,V}(G)\subset \mathbb{R}^2$ of all vectors $w$ such that  there exists a compact subset $K$ of the plane, and a sequence $(m_k,n_k,p_k)_{k\geq 0}$ of elements of $\mathbb{Z}^3$ so that:
\begin{enumerate}
\item for every $k$, $U^{-m_k}V^{-n_k}G^{p_k}(K) \cap K \neq \emptyset$,
\item the sequence $(m_k,n_k,p_k)$ tends to infinity,
\item the sequence $\left(\frac{m_k}{p_k},\frac{n_k}{p_k}\right)$ tends to $w$.
\end{enumerate}
\end{defi*}

\begin{rema}
In the case where $F$ is a lift of a homeomorphism of $\mathbb{T}^2$, and $S,T$ are the elementary translations $(x,y)\mapsto (x+1,y)$ and $(x,y)\mapsto (x,y+1)$, one easily checks that the rotation set $\rho_{S,T}(F)$ coincides with the classical rotation set of $F$.
\end{rema}

In order to prove the above theorem, we will consider a lift $F$ of a torus homeomophism whose rotation set (in the sense of Misiurewicz and Ziemiann) is the given compact convex set $C$. In order to realize the set $C'=L(C)$, we will not only replace $F$ by a new homeomophism $G$; we will also replace the elementary translations $S:(x,y)\mapsto (x+1,y)$ and $T:(x,y)\mapsto (x,y+1)$ by some ``non-linear translations" $U,V$.

\begin{rema}
The above definition immediatly extends to the case of a $\mathbb{Z}^p$ action on a non compact topological space $X$. In this more general setting, to get a more symmetric definition, it is tempting to replace $U^{-m_k}V^{-n_k}G^{p_k}$ in the first item by $U^{-m_k}V^{-n_k}G^{-p_k}$, and to define the ``rotation set" as a subset of $\mathbb{RP}^{p-1}$, instead of looking in a specific affine chart. The definition depends on a choice of basis of $\mathbb{Z}^p$, but two different choices give two ``rotation sets" that differ under a projective transformation, thus we get a conjugacy invariant which is a subset of  $\mathbb{RP}^{p-1}$ up to projective isomorphisms (see the argument at the end of the paper). Going back to the case of an action of $\mathbb{Z}^3$ on $\mathbb{R}^2$, one could wonder which results of the classical rotation set theory for torus homeomorphisms (in the sense of Misiurewicz and Ziemian) can be generalized to rotation sets of $\mathbb{Z}^3$ actions on the plane.
\end{rema}

Now we are in a position to give a more precise statement of the theorem above. We denote by $\Delta_\infty=\{[x:y:0]\}$ the ``line at infinity" in $\mathbb{R}\mathbb{P}^2$, and by $\Phi: \mathbb{R}\mathbb{P}^2\setminus \Delta_\infty\to \mathbb{R}^2$ the affine chart mapping $[x:y:z]$ to $(x/z,y/z)$. If $L\in \mathrm{SL}(3,\mathbb{R})$ is a projective transformation, we denote by $\widehat L$ the ``restriction of this map to the affine plane": more formally,
$$\widehat L=\Phi L \Phi^{-1} :\mathbb{R}^2\setminus (\Phi L^{-1}(\Delta_\infty)) \longrightarrow \mathbb{R}^2\setminus (\Phi L(\Delta_\infty)).$$

\begin{theo*}
Let $S : (x,y)\mapsto (x+1,y)$ and $T : (x,y)\mapsto (x,y+1)$. 
Let $F$ be a lift of a homeomorphism of the torus $\mathbb{T}^2 = \mathbb{R}^2/\mathbb{Z}^2$ isotopic to the identity.
Let $L \in \mathrm{SL}(3,\mathbb{Z})$ be a projective transformation such that $\rho_{S,T}(F)$ is disjoint from the line $\Phi L^{-1}(\Delta_\infty)$. Let 
$$L^{-1}=\left(\begin{array}{ccc}a_1 & a_2 & a_3 \\ b_1 & b_2 & b_3 \\ c_1 & c_2 & c_3 \end{array}\right)\quad\quad\mbox{and} \quad\quad \begin{array}{ccc} U = S^{a_1}T^{b_1}F^{-c_1}\\ V = S^{a_2}T^{b_2}F^{-c_2}\\ G= S^{-a_3}T^{-b_3}F^{c_3}\end{array}.\text{ Then:}$$ 
\begin{enumerate}
\item the quotient space $\mathbb{R}^2/\langle U,V\rangle$ is homeomorphic to the torus $\mathbb{T}^2$;
\item the rotation set $\rho_{U,V}(G)$ is equal to $\widehat L(\rho_{S,T}(F))$.
\end{enumerate}
\end{theo*}

\begin{rema}
Note that, since $G$ obviously commutes with $U$ and $V$, it can be seen as a lift of a homeomorphism $g$ of $\mathbb{T}^2$  which is isotopic to the identity. Thus this second theorem implies the first one.
From the definition, one easily deduces that $g$ is as smooth as $f$: if $f$ is $C^r$ for some $r\in \mathbb{N}\cup\{\infty\}$ or analytical, then so is $g$. Moreover, every invariant finite measure for $f$ induces an invariant finite measure for $g$. For example, if $f$ preserves a measure in the Lebesgue class, then so does $g$. 
\end{rema}

\begin{rema}
\label{r.surfaces-classification}
Consider a $\mathbb{Z}^2$ action on $\mathbb{R}^2$ generated by some homeomorphisms $U$ and $V$.  Assume this action is properly discontinuous. Then the quotient space $\mathbb{R}^2/\langle U,V\rangle$ is a topological surface (\emph{i.e.} a separated topological manifold of dimension 2)  whose fundamental group is isomorphic to $\mathbb{Z}^2$. According to the classification of surfaces (see \emph{e.g.} \cite{Ric}), it follows that this quotient space must be homeomorphic to $\mathbb{T}^2$. This is a key ingredient of the following proof that will play the part of Fried's theorem in Kwapisz's original proof.
\end{rema} 

\begin{proof}[Proof of Item 1 of the theorem.]
In view of Remark~\ref{r.surfaces-classification}, it is enough to prove that the action of $\mathbb{Z}^2$ on $\mathbb{R}^2$ generated by the homeomorphisms $U$ an $V$ is properly discontinuous: we consider a ball $B(0,R)$ in $\mathbb{R}^2$, and we aim to prove that $U^mV^n(B(0,R))$ is disjoint from $B(0,R)$ whenever $\|(m,n)\|$ is large enough.  

We denote by $\mathrm{D}(H)$ the \emph{displacement set} of the homeomorphism $H$ of the plane, that is, the set of all vectors of the type $H(z)-z$ where $z$ ranges over $\mathbb{R}^2$. Obviously $\mathrm{D}(S)=\{(1,0)\}$ and $\mathrm{D}(T)=\{(0,1)\}$. By assumption, the rotation set $\rho_{S,T}(F)$ is disjoint from the line $\Phi L^{-1}(\Delta_\infty)$. Therefore, we may consider a compact neighbourhood $\mathcal{O}$ of $\rho_{S,T}(F)$ so that 
$$\mathrm{dist}(\mathcal{O},\mathbb{R}^2\cap \Phi L^{-1}(\Delta_\infty))>\epsilon>0.$$
From the definition of the rotation set, one immediately sees that there exists an integer $k_0$ so that $\mathrm{D}(F^k)\subset k\mathcal{O}$ for $|k|\geq k_0$. And since $D(F^k)$ is bounded for every $|k|<k_0$, one gets that there exists $R'$ so that, for every $k\in\mathbb{Z}$ 
$$\mathrm{D}(F^k)\subset B(0,R')+k\mathcal{O}.$$
Now recall that $U = S^{a_1}T^{b_1}F^{-c_1}$ and $V = S^{a_2}T^{b_2}F^{-c_2}$. Since $S$, $T$ and $F$ commute, one immediately gets, for every $(m,n)\in\mathbb{Z}^2$,
$$\mathrm{D}(U^mV^n)=(ma_1+na_2,mb_1+nb_2)-\mathrm{D}(F^{mc_1+nc_2}).$$
Using the inclusion above, we obtain that, for every $(m,n)\in\mathbb{Z}^2$,
$$\mathrm{D}(U^mV^n)\subset B(0,R')+(ma_1+na_2,mb_1+nb_2)-(mc_1+nc_2) \mathcal{O},$$
and therefore 
\begin{eqnarray*}
U^mV^n(B(0,R)) & \subset & B(0,R+R')+(ma_1+na_2,mb_1+nb_2)-(mc_1+nc_2) \mathcal{O}\\
& = & B(0,R+R')+(mc_1+nc_2)\Phi L^{-1}([m:n:0])-(mc_1+nc_2)\mathcal{O}\\
& \subset & B(0,R+R')+(mc_1+nc_2)\left(\Phi L^{-1}(\Delta_\infty)-\mathcal{O}\right)\\
& \subset & B(0,R+R')+(mc_1+nc_2)(\mathbb{R}^2\setminus B(0,\epsilon)).
\end{eqnarray*}
(The last inclusion comes from the definition of the neighbourhood $\mathcal{O}$.)

On the first hand, if $|mc_1+nc_2|$ is  larger than $\frac{2R+R'}{\epsilon}$, the last inclusion above implies that $U^mV^n(B(0,R))$ is disjoint from $B(0,R)$, as desired. On the other hand, since $\mathcal{O}$ is compact, we can find $R''$ so that $(mc_1+nc_2) \mathcal{O}\subset B(0,R'')$ whenever $|mc_1+nc_2|\leq \frac{2R+R'}{\epsilon}$. As a consequence, if $|mc_1+nc_2|$ is smaller than $\frac{2R+R'}{\epsilon}$, but $\|(ma_1+na_2,mb_1+nb_2)\|$ is larger than $2R+R'+R''$, then the first inclusion above implies that $U^mV^n(B(0,R))$ is disjoint from $B(0,R)$ as desired.

To conclude, it remains to notice that since the vectors $(a_1,b_1,c_1)$ and $(a_2,b_2,c_2)$ are non-colinear (recall that $L^{-1}$ has rank three), the map $(m,n) \mapsto (ma_1+na_2,mb_1+nb_2,mc_1+nc_2)$ is a proper embedding of $\mathbb{Z}^2$ into $\mathbb{R}^3$. Thus the two quantities $|mc_1+nc_2|$ and $\|(ma_1+na_2,mb_1+nb_2)\|$ cannot remain bounded at the same time when $\|(m,n)\|$ is large. This shows that $U^mV^n(B(0,R))$ is disjoint from $B(0,R)$ provided that $\|(m,n)\|$ is bigger than some constant. In other words, the action of $\mathbb{Z}^2$ on $\mathbb{R}^2$ generated by the homeomorphisms $U$ an $V$ is properly discontinuous. According to Remark~\ref{r.surfaces-classification}, this implies that $\mathbb{R}^2/\langle U,V\rangle$ is homeomorphic to $\mathbb{T}^2$.
\end{proof}

\begin{proof}[Proof of Item 2 of the theorem.]
Consider a  compact subset $K$ of $\mathbb{R}^2$ and two sequences $(m_k,n_k,p_k)_{k\geq 0}$ and $(\mu_k,\nu_k,\pi_k)_{k\geq 0}$ of elements of $\mathbb{Z}^3$ which are related by $$(\mu_k,\nu_k,\pi_k)=L(m_k,n_k,p_k).$$
Obviously, $(m_k,n_k,p_k)$ tends to infinity if and only if $(\mu_k,\nu_k,\pi_k)$ tends to infinity. 
Now observe that $$S^{-m_k}T^{-n_k}F^{p_k}=U^{-\mu_k}V^{-\nu_k}G^{\pi_k}.$$ 
In particular, $(S^{-m_k}T^{-n_k}F^{p_k}(K))\cap K \neq \emptyset$ if and only if $(U^{-\mu_k}V^{-\nu_k}G^{\pi_k}(K))\cap K \neq \emptyset$.  Finally, $\left(m_k/p_k,n_k/p_k\right)$ converges to $w\in\mathbb{R}^2$ if and only if
$\left(\mu_k/\pi_k,\nu_k/\pi_k\right)=\widehat L\left(m_k/p_k,n_k/p_k\right)$ 
converges to the vector $\widehat L(w)$.
This shows that $\rho_{U,V}(G)=\widehat L(\rho_{S,T}(F))$.
\end{proof}

\bigskip

\noindent
\emph{Fran\c{c}ois B\'eguin}\\
LAGA, CNRS UMR 7539, Universit\'e Paris 13, 93430 Villetaneuse, France.
\smallskip

\noindent
\emph{Sylvain Crovisier}\\
LMO, CNRS UMR 8628, Universit\'e Paris-Sud 11, 91405 Orsay, France.
\smallskip

\noindent
\emph{Fr\'ed\'eric Le Roux}\\
IMJ-PRG, CNRS UMR 7586, Université Marie et Pierre Curie, 75005  Paris, France. 
\begin{thebibliography}{99}



\bibitem{FraMis}
Franks, John and Misiurewicz, Michal. Rotation sets of toral flows.  \textit{Proc. Amer. Math. Soc.} \textbf{109}  (1990),  no. 1, 243--249.

\bibitem{Kwa}
Kwapisz, Jaroslaw.
A  priori degeneracy  of  one-dimensional rotation  sets for  periodic
point free  torus maps. \textit{Trans. Amer.  Math. Soc.} \textbf{354}
(2002), no 7, 2865--2895.  


\bibitem{Kor}
Koropecki, Andres and Passeggi, Alejandro and Sambarino, Mart{\'\i}n.
The Franks-Misiurewicz conjecture for extensions of irrational rotations.
\textit{https://arxiv.org/abs/1611.05498}, 2016.


\bibitem{LecTal}
Le Calvez, Patrice and Tal, Fabio Armando.
Forcing theory for transverse trajectories of surface homeomorphisms.
\textit{https://arxiv.org/abs/1503.09127}, 2015.


\bibitem{MisZie}
Misiurewicz, Michal and Ziemian, Krystyna.
Rotation sets for maps of tori.
\textit{J. London. Math Soc.} (2) \textbf{40} (1989), no.3, 490--506.


\bibitem{Ric}
Ian Richards. 
On the classification of noncompact surfaces. 
\textit{Trans. Amer. Math. Soc.} \textbf{106} (1963), 259--269.


\end{thebibliography}
\end{document}